\documentclass{amsart}
\usepackage[latin1]{inputenc}
\usepackage[active]{srcltx} 
\newtheorem{thm}{Theorem}
\newtheorem{cor}[thm]{Corollary}
\newtheorem{lem}[thm]{Lemma}
\newtheorem{prop}[thm]{Proposition}
\theoremstyle{definition}
\newtheorem{defn}[thm]{Definition}
\theoremstyle{remark}

\numberwithin{equation}{section}

\newcommand{\To}{\longrightarrow}
\begin{document}

\thanks{Author supported by FPU grant of SEEU-MECD of Spain.}
\title[]{Radon-Nikod\'ym compact spaces of low weight and Banach spaces}%
\author{Antonio Avil\'es}
\address{Departamento de Matem\'aticas\\ Universidad de Murcia\\ 30100 Espinardo (Murcia)\\ Spain }%
\email{avileslo@um.es}

\begin{abstract}
We prove that a continuous image of a Radon-Nikodým compact of weight less
than $\mathfrak{b}$ is
Radon-Nikodým compact. As a Banach space counterpart,
subspaces of Asplund generated Banach spaces of density character
less than $\mathfrak{b}$ are Asplund generated. In this case, in addition,
there exists a subspace of an Asplund generated space which is not Asplund generated
which has density character exactly $\mathfrak{b}$.
\end{abstract}

\maketitle

The concept of Radon-Nikodým compact, due to Reynov~\cite{Reynov}, has its origin in Banach space
theory, and it is defined as a topological space which is homeomorphic to
a weak$^\ast$ compact subset of the dual of an Asplund space, that is, a
dual Banach space with the Radon-Nikodým property (topological spaces
will be here assumed to be Hausdorff). In~\cite{Namioka},
the following characterization of this class is given:\\

\begin{thm}\label{lsc metric}
A compact space $K$ is Radon-Nikodým compact if and only if there is a
lower semicontinuous metric $d$ on $K$ which fragments $K$.\\
\end{thm}

Recall that a map $f:X\times X\To \mathbb{R}$ on a topological
space $X$ is said to \emph{fragment} $X$ if for each (closed)
subset $L$ of $X$ and each $\varepsilon>0$ there is a nonempty
relative open subset $U$ of $L$ of $f$-diameter less than
$\varepsilon$, i.e. $\sup\{f(x,y) : x,y\in U\}<\varepsilon$. Also,
a map $g:Y\To\mathbb{R}$ from a topological space to the real line
is \emph{lower semicontinuous} if
$\{y : g(y)\leq r\}$ is closed in $Y$ for every real number $r$.\\

It is an open problem whether a continuous image of a Radon-Nikod\'ym
compact is Radon-Nikod\'ym. Arvanitakis~\cite{Arvanitakis} has made the
following approach to this problem:
if $K$ is a Radon-Nikodým compact and $\pi:K\To L$ is a continuous surjection, then we have
a lower semicontinuous
fragmenting metric $d$ on $K$, and if we want to prove that $L$ is
Radon-Nikodým compact, we should find such a metric on $L$. A natural candidate
is:
\begin{eqnarray*}
d_1(x,y) &=& d(\pi^{-1}(x),\pi^{-1}(y)) = \inf\{d(t,s) : \pi(t)=x,\
\pi(s)=y\}.
\end{eqnarray*}

The map $d_1$ is lower semicontinuous and fragments $L$ and it is
a \emph{quasi metric}, that is, it is symmetric and vanishes only
if $x=y$. But it is not a metric because, in general, it lacks
triangle inequality. Consequently, Arvanitakis~\cite{Arvanitakis}
introduced the following concept:\\

\begin{defn}
A compact space $L$ is said to be \emph{quasi Radon-Nikodým} if
there exists a lower semicontinuous
quasi metric which fragments $L$.\\
\end{defn}

The class of quasi Radon-Nikodým compacta is closed under
continuous images but it is unknown whether it is the same class
as that of Radon-Nikodým compacta or even the class of their
continuous images. At least other two superclasses of continuous
images of Radon-Nikodým compacta appear in the literature.
Reznichenko~\cite[p. 104]{ArkhangelskiiGT2} defined a compact
space $L$ to be \emph{strongly fragmentable} if there is a metric
$d$ which fragments $L$ such that each pair of different points of
$L$ possess disjoint neighbourhoods at a positive $d$-distance. It
has been noted by Namioka~\cite{Namiokanote} that the classes of
quasi Radon-Nikodým and strongly fragmentable compacta are equal.
The other superclass of continuous images of Radon-Nikodým
compacta, called \emph{countably lower fragmentable} compacta, was
introduced by Fabian, Heisler and Matou\v{s}ková~\cite{FabHeiMat}.
In section~\ref{clf}, we recall its definition and we prove that
this class is equal to the other two.\\

The main result in section~\ref{cubos} is the following:\\

\begin{thm}\label{cardinal b}
If $K$ is a quasi Radon-Nikodým compact space of weight less than
$\mathfrak{b}$, then $K$ is Radon-Nikodým compact.\\
\end{thm}

The weight of a topological space is the least cardinality of a
base for its topology. We also recall the definition of cardinal
$\mathfrak{b}$. In the set $\mathbb{N}^\mathbb{N}$ we consider the
order relation given by $\sigma\leq\tau$ if $\sigma_n\leq \tau_n$
for all $n\in\mathbb{N}$. Cardinal $\mathfrak{b}$ is the least
cardinality of a subset of $\mathbb{N}^\mathbb{N}$ which is not
\emph{$\sigma$-bounded} for this order (a set is
\emph{$\sigma$-bounded} if it is a countable union of bounded
subsets). It is consistent that $\mathfrak{b}>\omega_1$. In fact,
Martin's axiom and the negation of the continuum hypothesis imply
that $\mathfrak{c}=\mathfrak{b}>\omega_1$, cf.~\cite[11D and
14B]{Martinsaxiom}. It is also possible that
$\mathfrak{c}>\mathfrak{b}>\omega_1$, cf.~\cite[section 5]{vanDouwen}.
On the other hand, cardinal $\mathfrak{d}$ is the least cardinality
of a cofinal subset of $(\mathbb{N}^\mathbb{N},\leq)$, that is, a
set $A$ such that for each $\sigma\in\mathbb{N}^\mathbb{N}$ there
is some $\tau\in A$ such that $\sigma\leq \tau$. In a sense, the
following proposition puts a rough bound on the size of the class
of quasi Radon-Nikodým compacta with respect to
Radon-Nikodým compacta.\\

\begin{prop}\label{cardinal d}
Every quasi Radon-Nikodým compact space embeds into a product of
Radon-Nikodým compact spaces with at most $\mathfrak{d}$ factors.\\
\end{prop}

In section~\ref{Banach} we discuss the Banach space counterpart to Theorem~\ref{cardinal b}.
A Banach space $V$ is Asplund generated, or $GSG$, if
there is some Asplund space $V'$ and a bounded linear operator $T:V'\To V$
such that $T(V')$ is dense in $V$. Our main result for this class is the
following:\\

\begin{thm}\label{Asplund generated}
Let $V$ be a Banach space of density character
less than $\mathfrak{b}$ and such that the dual unit ball
$(B_{V^\ast},w^\ast)$ is quasi Radon-Nikod\'ym compact, then $V$ is Asplund generated.\\
\end{thm}

The density character of a Banach space is the least cardinal of a
norm-dense subset, and it equals the weight of its dual unit ball in the
weak$^\ast$ topology.\\

Examples constructed by
Rosenthal~\cite{SubspacesWCG} and Argyros~\cite[section~1.6]{FabianWA}
show that there exist Banach spaces which are
subspaces of Asplund generated spaces but which are not Asplund generated.
However, since
the dual unit ball of a subspace of an Asplund generated space is a continuous image of a
Radon-Nikod\'ym compact \cite[Theorem~1.5.6]{FabianWA}, we have the
following corollary to Theorem~\ref{Asplund generated}:\\

\begin{cor}
If a Banach space $V$ is a subspace of an Asplund generated space and the density character of $V$
is less than $\mathfrak{b}$, then $V$ is Asplund generated.\\
\end{cor}

Also, a Banach space is weakly compactly generated (WCG) if it is the
closed linear span of a weakly compact subset.
The same examples mentioned above show that neither is this property
inherited by subspaces. A Banach space $V$ is weakly compactly generated if and only if it is
Asplund generated and its dual unit ball $(B_{V^\ast},w^\ast)$ is Corson
compact~\cite{OriSchVal}, \cite{StegallRN}.
Having Corson dual unit ball is a hereditary property
since a continuous image of a Corson compact is Corson
compact~\cite{CorsonImages}, hence:\\

\begin{cor}\label{subWCG}
If a Banach space $V$ is a subspace of a weakly compactly generated space and the density character of $V$
is less than $\mathfrak{b}$, then $V$ is weakly compactly generated.\\
\end{cor}

Corollary~\ref{subWCG} can also be obtained from the following theorem,
essentially due to Mercourakis~\cite{MerSta}:\\

\begin{thm}\label{WKA below b}
If a Banach space $V$ is weakly $\mathcal{K}$-analytic and the density character of $V$
is less than $\mathfrak{b}$, then $V$ is weakly compactly generated.\\
\end{thm}

The class of weakly $\mathcal{K}$-analytic spaces is larger than the class of subspaces
of weakly compactly generated spaces. We recall its definition in section~\ref{Banach}.
The result of Mercourakis~\cite[Theorem 3.13]{MerSta} is that, under Martin's axiom,
weakly $\mathcal{K}$-analytic Banach spaces of density character less than $\mathfrak{c}$
are weakly compactly generated, but his arguments prove in fact the more
general Theorem~\ref{WKA below b}. We give a more elementary proof of this
theorem, obtaining it as a consequence of a purely topological result: Any
$\mathcal{K}$-analytic topological space of density character less than
$\mathfrak{b}$ contains a dense $\sigma$-compact subset. We also remark that it is not
possible to generalize Theorem~\ref{WKA below b} for the class of weakly countably
determined Banach spaces.\\

Cardinal $\mathfrak{b}$ is best possible for
Theorem~\ref{Asplund generated}, Theorem~\ref{WKA below b}
and their corollaries, as it is shown by slight modifications of
the
mentioned example of Argyros~\cite[section 1.6]{FabianWA}
and of
 the
example of Talagrand~\cite{WKA} of a weakly $\mathcal{K}$-analytic Banach
space which is not weakly compactly generated,
so that we get
examples of density character exactly $\mathfrak{b}$.\\

For information about cardinals $\mathfrak{b}$ and $\mathfrak{d}$ we refer
to~\cite{vanDouwen}. Concerning Banach spaces, our main reference
is~\cite{FabianWA}.\\

I want to express my gratitude to Jos\'e Orihuela for valuable
discussions and suggestions and to Witold Marciszewski, from whom
I learnt about cardinals $\mathfrak{b}$ and $\mathfrak{d}$. I also
thank Isaac Namioka and the referee for suggestions which have
improved the final version of this article.\\

\section{Quasi Radon-Nikodým compacta of low weight}\label{cubos}

In this section, we characterize quasi Radon-Nikodým compacta in
terms of embeddings into cubes $[0,1]^\Gamma$ and from this, we
will derive proofs of Theorem~\ref{cardinal b} and
Proposition~\ref{cardinal d}. Techniques of
Arvanitakis~\cite{Arvanitakis} will play an important role in this
section, as well as the following theorem of Namioka~\cite{Namioka}:\\

\begin{thm}
Let $K$ be a compact space. The following are equivalent.
\begin{enumerate}
\item $K$ is Radon-Nikodým compact.
\item There is an embedding $K\subset [0,1]^\Gamma$ such that $K$ is
fragmented by the uniform metric $d(x,y) = \sup_{\gamma\in\Gamma}|x_\gamma
- y_\gamma|$.\\
\end{enumerate}
\end{thm}

Let $\mathrm{P}\subset\mathbb{N}^\mathbb{N}$ be the set of all
strictly increasing sequences of positive integers. Note that this
is a cofinal subset of $\mathbb{N}^\mathbb{N}$. For each
$\sigma\in\mathrm{P}$ we consider the lower semicontinuous
non decreasing function $h^\sigma:[0,+\infty]\To \mathbb{R}$ given by:\\

\begin{itemize}
\item $h^\sigma(0) = 0$,\\

\item $h^\sigma(t) = \frac{1}{\sigma_n}$ whenever $\frac{1}{n+1}<t\leq\frac{1}{n}$.\\

\item $h^\sigma(t)=\frac{1}{\sigma_1}$ whenever $t>1$.\\
\end{itemize}

Also, if $f:X\times X\To \mathbb{R}$ is a map and $A,B\subset X$, we
will use the notation $f(A,B) = \inf\{f(x,y) : x\in A,\ y\in B\}$.\\

\begin{thm}\label{qRN embedding}
Let $K$ be a compact subset of the cube $[0,1]^\Gamma$. The
following are equivalent:
\begin{enumerate}
\item $K$ is quasi Radon-Nikodým compact.

\item There is a map $\sigma:\Gamma\To\mathrm{P}$ such that $K$ is
fragmented by
$$f(x,y) = \sup_{\gamma\in\Gamma}h^{\sigma(\gamma)}(|x_\gamma -
y_\gamma|)$$ which is a lower semicontinuous quasi metric.\\
\end{enumerate}
\end{thm}

PROOF: Observe that $f$ in (2) is expressed as a supremum of lower
semicontinuous functions, and therefore, it is lower
semicontinuous. Also, $f(x,y)=0$ if and only if
$h^{\sigma(\gamma)}(|x_\gamma - y_\gamma|)=0$ for all
$\gamma\in\Gamma$ if and only if $|x_\gamma-y_\gamma|=0$ for all
$\gamma\in\Gamma$. Hence, $f$ is indeed a lower semicontinuous
quasi metric and it is clear that (2) implies (1). Assume now that
$K$ is quasi Radon-Nikodým compact and let $g:K\times K\To [0,1]$
be a lower semicontinuous quasi metric which fragments $K$. For
$\gamma\in\Gamma$, we call $p_\gamma:K\To [0,1]$ the projection on
the coordinate $\gamma$, $p_\gamma(x)=x_\gamma$, and we define a
quasi metric $g_\gamma$ on $[0,1]$ by the rule:\\

$g_\gamma (t,s)=
  \begin{cases}
    g(p_\gamma^{-1}(t),p_\gamma^{-1}(s)) & \text{if }p_\gamma^{-1}(t)\text{ and }p_\gamma^{-1}(s)\text{ are nonempty}, \\
    1 & \text{otherwise}.
  \end{cases}$

\ \\
Note that $g_\gamma$ is lower semicontinuous because for $r<1$

$$ \{(t,s) : g_\gamma(t,s)\leq r\} = \bigcap_{r'> r}(p_\gamma\times p_\gamma)\{(x,y)\in K^2 : g(x,y)\leq r'\}$$

Observe also that if $x,y\in K$, then $g_\gamma(x_\gamma,y_\gamma)= g_\gamma(p_\gamma(x),p_\gamma(y))\leq
g(x,y)$. Hence, $K$ is fragmented by
$$ g'(x,y) = \sup_{\gamma\in\Gamma}g_\gamma(x_\gamma,y_\gamma) \leq
g(x,y)$$

The proof finishes by making use of the following lemma, where we put $g_0:=g_\gamma$:\\

\begin{lem}
Let $g_0:[0,1]\times[0,1]\To[0,1]$ be a lower semicontinuous quasi metric on $[0,1]$. Then, there exists
$\tau\in\mathrm{P}$ such that $h^\tau(|t-s|)\leq g_0(t,s)$ for all
$t,s\in [0,1]$.
\end{lem}

PROOF: We define $\tau$ recursively. Suppose that we have defined
$\tau_1,\ldots,\tau_n$ in such a way that if $|t-s|>\frac{1}{n+1}$,
then $h^\tau(|t-s|)\leq g_0(t,s)$. Let
$$K_m = \left\{(t,s)\in [0,1]\times [0,1] : |t-s|\geq\frac{1}{n+2} \text{ and }
 g_0(t,s)\leq\frac{1}{m}\right\}$$

Then, $\{K_m\}_{m=1}^\infty$ is a decreasing sequence of compact subsets of
$[0,1]^2$
with empty intersection. Hence, there is $m_1$ such that $K_{m}$ is empty for $m\geq m_1$.
We define $\tau_{n+1} = \max\{m_1, \tau_n+1\}$.\hfill$\square$\\

Now,
we state a lemma which
is just a piece of the proof of~\cite[Proposition 3.2]{Arvanitakis}.
We include its proof for the sake of completeness.\\

\begin{lem}\label{imagen continua}
Let $K,L$ be compact spaces, let $f:K\times K\To\mathbb{R}$ be a
symmetric map which fragments $K$ and $p:K\To L$ a continuous
surjection. Then $L$ is fragmented by $g(x,y) =
f(p^{-1}(x),p^{-1}(y))$ and in particular, $L$ is fragmented by
any $g'$ with $g'\leq g$.
\end{lem}

PROOF: Let $M$ be a closed subset of $L$ and $\varepsilon>0$. By Zorn's
lemma a set $N\subset K$ can be found such that $p:N\To M$ is onto and
irreducible (that is, for every $N'\subset N$ closed, $p:N'\To M$ is not
onto). We find $U\subset N$ a relative open subset of $N$ of $f$-diameter less than
$\varepsilon$. By irreducibility, $p(U)$ has nonempty relative interior in $M$.
This interior is a nonempty relative open subset of $M$ of $g$-diameter less than
$\varepsilon$.\hfill$\square$\\

In the sequel, we use the following notation: If $A\subset\Gamma$ are
sets, $d_A$ states for the pseudometric in $[0,1]^\Gamma$ given by
$d_A(x,y) = \sup_{\gamma\in A}|x_\gamma-y_\gamma|$.\\

\begin{lem}\label{bounded set}
Let $K$ be a compact subset of the cube $[0,1]^\Gamma$ and let
$\sigma:\Gamma\To\mathrm{P}$ be a map such that the quasi metric
$$f(x,y) = \sup_{\gamma\in\Gamma}h^{\sigma(\gamma)}(|x_\gamma -
y_\gamma|)$$ fragments $K$ and such that $\sigma(\Gamma)$ is a
$\sigma$-bounded subset of $\mathbb{N}^\mathbb{N}$. Then, $K$ is
Radon-Nikodým compact. In addition, there exist sets
$\Gamma_n\subset\Gamma$ such that
$\Gamma=\bigcup_{n\in\mathbb{N}}\Gamma_n$
and each $d_{\Gamma_n}$ fragments $K$.\\
\end{lem}

PROOF: There is a decomposition $\Gamma =
\bigcup_{n\in\mathbb{N}}\Gamma_n$ such that each
$\sigma(\Gamma_n)$ has a bound $\tau_n$ in
$(\mathbb{N}^\mathbb{N},\leq)$. We choose $\tau_n\in\mathrm{P}$.
First, we prove that each $d_{\Gamma_n}$ fragments $K$. For every
$n\in\mathbb{N}$, $K$ is fragmented by the map
$$f_n(x,y) = \sup_{\gamma\in\Gamma_n}h^{\sigma(\gamma)}(|x_\gamma -
y_\gamma|)\leq f(x,y)$$
 and
\begin{eqnarray*} f_n(x,y) &=&
\sup_{\gamma\in\Gamma_n}h^{\sigma(\gamma)}(|x_\gamma - y_\gamma|)
\geq \sup_{\gamma\in\Gamma_n}h^{\tau_n}(|x_\gamma - y_\gamma|)\\
&=& h^{\tau_n}\left(\sup_{\gamma\in\Gamma_n}|x_\gamma - y_\gamma|\right)
 =h^{\tau_n}(d_{\Gamma_n}(x,y)).\\
\end{eqnarray*}

Hence, a set of
$f_n$-diameter less than $\frac{1}{\tau_n}$ in $K$ is a set of $d_{\Gamma_n}$-diameter
less than
$\frac{1}{n}$ and therefore, since $f_n$ fragments $K$, also $d_{\Gamma_n}$ fragments $K$.\\

Consider now $p_n:[0,1]^\Gamma\To [0,1]^{\Gamma_n}$ the natural projection and
$K_n = p_n(K)$. By Lemma~\ref{imagen continua}, since $K$ is fragmented by $f_n$, $K_n$ is fragmented
by $$g_n(x,y) = \sup_{\gamma\in\Gamma_n}h^{\sigma(\gamma)}(|x_\gamma -
y_\gamma|).$$ and hence, $K_n$ is Radon-Nikod\'ym compact. Moreover, since
$\Gamma=\bigcup_{n\in\mathbb{N}}\Gamma_n$, $K$ embeds
into the product $\prod_{n\in\mathbb{N}}K_n$ and the class of
Radon-Nikodým compacta is closed under taking countable products and under
taking closed subspaces~\cite{Namioka}, so $K$ is Radon-Nikod\'ym
compact.\qed\\

PROOF OF THEOREM~\ref{cardinal b}: If the weight of $K$ is less
than $\mathfrak{b}$, then $K$ can be embedded into a cube
$[0,1]^\Gamma$ with $|\Gamma|<\mathfrak{b}$. Any subset of
$\mathbb{N}^\mathbb{N}$ of cardinality less than $\mathfrak{b}$ is
$\sigma$-bounded, so Theorem~\ref{cardinal b} follows directly
from
Theorem~\ref{qRN embedding} and Lemma~\ref{bounded set}.\hfill$\square$\\

PROOF OF PROPOSITION~\ref{cardinal d}: Let $K$ be quasi
Radon-Nikodým compact, suppose $K$ is embedded into some cube
$[0,1]^\Gamma$ and let $\sigma:\Gamma\To\mathrm{P}$ be as in
Theorem~\ref{qRN embedding}. Let $A\subset\mathrm{P}$ be a cofinal
subset of $\mathrm{P}$ of cardinality $\mathfrak{d}$. For $\alpha\in
A$, let
$$\Gamma_\alpha = \{\gamma\in\Gamma : \sigma(\gamma)\leq\alpha\},$$
let $p_\alpha:[0,1]^\Gamma\To [0,1]^{\Gamma_\alpha}$ be the natural
projection, and let $K_\alpha = p_\alpha(K)$.
Again, since $\Gamma = \bigcup_{\alpha\in A}\Gamma_\alpha$, $K$ embeds
into the product $\prod_{\alpha\in A}K_\alpha$. By Lemma~\ref{imagen
continua}, $K_\alpha$ is fragmented by
$$g_\alpha(x,y) =
\sup_{\gamma\in\Gamma_\alpha}h^{\sigma(\gamma)}(|x_\gamma-y_\gamma|)$$
The set $\{\sigma(\gamma) :\gamma\in\Gamma_\alpha\}$ is a bounded,
and hence $\sigma$-bounded, set. Hence, by Lemma~\ref{bounded
set},
$K_\alpha$ is Radon-Nikodým compact.\hfill $\square$\\

We note that from Lemma~\ref{bounded set}, we obtain something
stronger than Theorem~\ref{cardinal b}:\\

\begin{thm}\label{countable pieces} For every quasi Radon-Nikod\'ym compact subset of a cube
$[0,1]^\Gamma$ with $|\Gamma|<\mathfrak{b}$ there is a countable
decomposition $\Gamma=\bigcup_{n\in\mathbb{N}}\Gamma_n$ such that
$d_{\Gamma_n}$ fragments $K$ for all $n\in\mathbb{N}$.\\ \end{thm}

A similar result holds also for generalized Cantor cubes
(cf.~\cite[Theorem 3]{FabHeiMat}, \cite[Theorem
3.6]{Arvanitakis}): If $K$ is a quasi Radon-Nikodým compact subset
of $\{0,1\}^\Gamma$, then there is a decomposition $\Gamma =
\bigcup_{n\in\mathbb{N}}\Gamma_n$ such that $d_{\Gamma_n}$
fragments $K$ for all $n\in\mathbb{N}$. We give now an example
which shows that this phenomenon does not happen for general
cubes, even if the compact $K$ has weight
less than $\mathfrak{b}$ or it is zero-dimensional:\\

\begin{prop}\label{counterexample pieces}
There exist a set $\Gamma$ of cardinality $\mathfrak{b}$ and a
compact subset $K$ of $[0,1]^\Gamma$ homeomorphic to the
metrizable Cantor cube $\{0,1\}^\mathbb{N}$ such that for any
decomposition $\Gamma=\bigcup_{n\in\mathbb{N}}\Gamma_n$ there
exists
$n\in\mathbb{N}$ such that $d_{\Gamma_n}$ does not fragment $K$.\\
\end{prop}

PROOF: First, we take $\Gamma$ a subset of $\mathbb{N}^\mathbb{N}$
of cardinality $\mathfrak{b}$ which is not $\sigma$-bounded. We call
$A = \{\gamma_n : \gamma\in\Gamma, n\in\mathbb{N}\}$ the set of
all terms of elements of $\Gamma$. We define
$$ K' = \{ x\in\{0,1\}^{\Gamma\times\mathbb{N}} :
x_{\gamma,n}=x_{\gamma',n'}\text{ whenever }\gamma_n =
\gamma'_{n'}\}.$$

Observe that $K'$ is homeomorphic to $\{0,1\}^\mathbb{N}$: namely,
for each $a\in A$ choose some
$\gamma^a,n^a\in\Gamma\times\mathbb{N}$ such that $\gamma^a_{n^a}=
a$; in this case we have a homeomorphism $K'\To \{0,1\}^A$ given
by $x\mapsto (x_{\gamma^a,n^a})_{a\in A}$.

Now, we consider the embedding
$\phi:\{0,1\}^{\Gamma\times\mathbb{N}}\To [0,1]^\Gamma$ given by
$$\phi(x) = \left(\sum_{n\in\mathbb{N}}\left(\frac{2}{3}\right)^n x_{\gamma,n}\right)_{\gamma\in\Gamma}$$

We claim that the space $K = \phi(K')\subset[0,1]^\Gamma$ verifies
the statement. Let $\Gamma = \bigcup_{n\in\mathbb{N}}\Gamma_n$ any
countable decomposition of $\Gamma$. Since $\Gamma$ is not
$\sigma$-bounded, there is some $n\in\mathbb{N}$ such that
$\Gamma_n$ is not bounded. For this fixed $n$, since $\Gamma_n$ is
not bounded, there is some $m\in\mathbb{N}$ such that the set $S =
\{\gamma_m : \gamma\in\Gamma_n\}\subset A$ is infinite. We
consider
$$ K_0 = \{x\in K' : x_{\gamma,k} = 0\text{ whenever }
\gamma_k\not\in S\}\subset K.$$ By the same arguments as for $K'$,
$K_0$ is homeomorphic to the Cantor cube $\{0,1\}^\mathbb{N}$ by a
map $K_0\To \{0,1\}^S$ given by $x\mapsto (x_{\gamma^a,n^a})_{a\in
S}$. Now, we take two different elements $x,y\in K_0$. Then, there
must exist some $\gamma\in\Gamma_n$ such that $x_{\gamma,m}\neq
y_{\gamma,m}$, and this implies that $|\phi(x)_\gamma -
\phi(y)_\gamma|\geq 3^{-m}$ and therefore
$d_{\Gamma_n}(\phi(x),\phi(y))\geq 3^{-m}$. This means that any
nonempty subset of $\phi(K_0)$ of $d_{\Gamma_n}$-diameter less
than $3^{-m}$ must be a singleton. If $d_{\Gamma_n}$ fragmented
$K$, this would imply that $\phi(K_0)$ has an isolated point,
which contradicts the fact that
it is homeomorphic to $\{0,1\}^\mathbb{N}$.$\qed$\\

\section{Banach spaces of low density character}\label{Banach}

In this section we find that cardinal $\mathfrak{b}$ is the least
possible density character of Banach spaces which are
counterexamples to several questions. First, we introduce some
notation: If $A$ is a subset of a Banach space $V$, we call $d_A$
to the pseudometric $d_A(x^\ast,y^\ast) = \sup_{x\in
A}|x^\ast(x)-y^\ast(x)|$ on
$B_{V^\ast}$. Also, we recall the following definition~\cite[Definition 1.4.1]{FabianWA}:\\

\begin{defn}
A nonempty bounded subset $M$ of a Banach space $V$ is called an
\emph{Asplund set} if for each countable set $A\subset M$ the
pseudometric space
$(B_{V^\ast},d_A)$ is separable.\\
\end{defn}

By~\cite[Theorem 2.1]{CasNamOri}, $M$ is an Asplund subset of $V$ if and only if $d_M$ fragments
$(B_{V^*},w^\ast)$. Also, by \cite[Theorem 1.4.4]{FabianWA}, a Banach space $V$ is Asplund
generated if and only if it is the closed linear span of an Asplund
subset.\\

PROOF OF THEOREM~\ref{Asplund generated}: Let $\Gamma$ be a dense
subset of the unit ball $B_V$ of $V$ of cardinality less than
$\mathfrak{b}$. Then, we have a natural embedding
$(B_{V^\ast},w^\ast)\subset [-1,1]^\Gamma$. Since
$(B_{V^\ast},w^\ast)$ is quasi Radon-Nikod\'ym compact, we apply
Theorem~\ref{countable pieces} and we have $\Gamma=\bigcup
\Gamma_n$ and each $d_{\Gamma_n}$ fragments $(B_{V^\ast},w^\ast)$.
This means that for each $n$, $\Gamma_n$ is an Asplund set, and
by~\cite[Lemma 1.4.3]{FabianWA}, $M =
\bigcup_{n\in\mathbb{N}}\frac{1}{n}\Gamma_n$ is an Asplund set
too. Finally, since the closed linear span of $M$ is $V$,
by~\cite[Theorem 1.4.4]{FabianWA}, $V$ is Asplund generated.\hfill $\square$\\

We recall now the concepts that we need for the proof of Theorem~\ref{WKA
below b}. We follow the terminology and notation
of~\cite[sections 3.1, 4.1]{FabianWA}. Let $X$ and $Y$ be topological spaces.
A map $\phi:X\to 2^Y$
from $X$ to the subsets of $Y$ is said to be an usco if the following conditions hold:
\begin{enumerate}
\item $\phi(x)$ is a compact
subset of $Y$ for all $x\in X$.
\item $\{x : \phi(x)\subset U\}$ is open in $X$, for every open set $U$ of
$Y$.
\end{enumerate}
In this situation, for $A\subset X$ we denote $\phi(A)=\bigcup_{x\in
A}\phi(x)$.\\

A completely regular topological space $X$ is said to be
$\mathcal{K}$-analytic if there exists an usco
$\phi:\mathbb{N}^\mathbb{N}\to 2^X$ such that
$\phi(\mathbb{N}^\mathbb{N}) = X$. A Banach space is weakly
$\mathcal{K}$-analytic if it is a
$\mathcal{K}$-analytic space in its weak topology.\\

We note that if a Banach space $V$ contains a weakly $\sigma$-compact
subset $M$ which is dense in the weak topology, then $V$ is WCG. This is
because if $M=\bigcup_{n=1}^\infty K_n$ being $K_n$ a weakly compact set
bounded by $c_n>0$, then $\{0\}\cup\bigcup\frac{1}{nc_n}K_n$ is a weakly
compact subset of $V$ whose linear span is (weakly) dense in $V$.
Hence, Theorem~\ref{WKA below b} is deduced from the
following:\\

\begin{prop}
If $X$ is a $\mathcal{K}$-analytic topological space which
contains a dense subset of cardinality less than $\mathfrak{b}$,
then $X$ contains a dense
$\sigma$-compact subset.\\
\end{prop}

PROOF: We have an usco $\phi:\mathbb{N}^\mathbb{N}\To 2^X$ with
$\phi(\mathbb{N}^\mathbb{N}) = X$ and also a set
$\Sigma\subset\mathbb{N}^\mathbb{N}$ such that
$|\Sigma|<\mathfrak{b}$ and $\phi(\Sigma)$ is dense in $X$. Any
subset of $\mathbb{N}^\mathbb{N}$ of cardinal less than
$\mathfrak{b}$ is contained in a $\sigma$-compact subset of
$\mathbb{N}^\mathbb{N}$~\cite[Theorem 9.1]{vanDouwen}. Uscos send
compact sets onto compact sets, so if $\Sigma'\supset\Sigma$ is
$\sigma$-compact,
then $\phi(\Sigma')$ is a dense $\sigma$-compact subset of $X$.\qed\\

We recall that a completely regular topological space $X$ is
$\mathcal{K}$-countably determined if there exists a subset
$\Sigma$ of $\mathbb{N}^\mathbb{N}$ and an usco $\phi:\Sigma\To
2^X$ such that $\phi(\Sigma)=X$ and that a Banach space is weakly
countably determined if it is $\mathcal{K}$-countably determined
in its weak topology. Talagrand~\cite{WCD} has constructed a
Banach space which is weakly countably determined but which is not
weakly $\mathcal{K}$-analytic. A slight modification of this
example gives a similar one with density character $\omega_1$.
This shows that no analogue of Theorem~\ref{WKA below b} is
possible for weakly countably determined Banach spaces. The change
in the example consists in substituting the set $T$ considered
in~\cite[p. 78]{WCD} by any subset $T'\subset T$ of cardinal
$\omega_1$ such that $\{o(X) : X\in T'\}$ is uncountable and
$\mathcal{A}$ by $\mathcal{A'}=\{A\subset T' :
A\in\mathcal{A}_1\}$ (the notations are explained in~\cite{WCD}).\\

Now, we turn to the fact that cardinal $\mathfrak{b}$ is best possible in
Theorem~\ref{Asplund generated}, Theorem~\ref{WKA below b} and their corollaries.
We fix a subset $S$ of $\mathbb{N}^\mathbb{N}$ of cardinality $\mathfrak{b}$ which is
not $\sigma$-bounded.\\

 Following
the exposition of the example of Argyros in~\cite[section 1.6]{FabianWA}
we just substitute the space $Y=\overline{span}\{\pi_\sigma :
\sigma\in\mathbb{N}^\mathbb{N}\}$ in~\cite[Theorem 1.6.3]{FabianWA}  by
$Y'=\overline{span}\{\pi_\sigma : \sigma\in S\}$ and
we obtain a Banach space of density character $\mathfrak{b}$ which is a
subspace of a WCG space $C(K)$ but which is not Asplund generated.
The same arguments in~\cite[section 1.6]{FabianWA} hold
just changing $\mathbb{N}^\mathbb{N}$ by $S$ where necessary. Only the proof
of~\cite[Lemma 1.6.1]{FabianWA} is not good for this case. It must be
substituted by the following:\\

\begin{lem}\label{new161}
Let $\Gamma_n$, $n\in\mathbb{N}$, be any subsets of $S$ such that
$\bigcup_{n\in\mathbb{N}}\Gamma_n = S$. Then there exist
$n,m\in\mathbb{N}$ and an
infinite set $A\in\mathcal{A}_m$ such that $A\subset\Gamma_n$.\\
\end{lem}

Here, as in~\cite[section 1.6]{FabianWA}, $\mathcal{A}_m$ is the family of
all subsets $A\subset\mathbb{N}^\mathbb{N}$ such that if $\sigma,\tau\in A$ and
$\sigma\neq\tau$, then $\sigma_i=\tau_i$ if $i\leq m$ and $\sigma_{m+1}\neq\tau_{m+1}$.
Also, $\mathcal{A}=\bigcup_{m=1}^\infty\mathcal{A}_m$.\\

PROOF OF LEMMA~\ref{new161}: We consider $\Gamma_{i,j} =
\{\sigma\in\Gamma_i : \sigma_1=j\}$, $i,j\in\mathbb{N}$. Note that
$S = \bigcup_{i,j}\Gamma_{i,j}$. Since $S$ is not
$\sigma$-bounded, there exist $n,l$ with $\Gamma_{n,l}$ unbounded.
This implies that for some $m$, the set $\{\sigma_m :
\sigma\in\Gamma_{n,l}\}$ is infinite. We take $m$ the least
integer with this property ($m>1$). Let $B\subset\Gamma_{n,l}$ be
an infinite set such that $\sigma_m\neq\sigma'_m$ for
$\sigma,\sigma'\in B$, $\sigma\neq\sigma'$. Since all $\sigma_k$
with $\sigma\in B$, $k<m$, lie in a finite set, an infinite set
$A\subset B$ can be
chose such that $A\in\mathcal{A}_{m-1}$.$\qed$\\

On the other hand, if we follow the proof in~\cite[section 4.3]{FabianWA}
that the Banach space $C(K)$ of Talagrand is weakly $\mathcal{K}$-analytic but not WCG,
and we change $K$ in~\cite[p.~76]{FabianWA} by $K' = \{\chi_{A} :
A\in\mathcal{A},\ A\subset S\}\subset\{0,1\}^S$ then $C(K')$ still
verifies this conditions and has density character $\mathfrak{b}$. Observe
that $C(K')$ is weakly $\mathcal{K}$-analytic because $K'$ is a retract of
the original $K$. The fact that $C(K')$ is not WCG (not
even a subspace of a WCG space) follows
from~\cite[Theorem~4.3.2]{FabianWA} and Lemma~\ref{new161} above by the
same arguments as in~\cite[p. 78]{FabianWA}.\\

\section{Countably lower fragmentable compacta}\label{clf}

In this section we prove that the concept of quasi Radon Nikodým
compact~\cite{Arvanitakis} is equivalent to that of countably
lower fragmentable compact~\cite{FabHeiMat}. The main result for
this class in~\cite{FabHeiMat} is that if $K$ is countably lower
fragmentable, then so is $(B_{C(K)^\ast},w^\ast)$. We note that,
with these two facts at hand, together with the fact that if
$C(K)$ is Asplund generated, then $K$ is
Radon-Nikod\'ym~\cite[Theorem~1.5.4]{FabianWA},
Theorem~\ref{cardinal b}
is deduced from Theorem~\ref{Asplund generated}.\\

We need some notation: if $K$ is a compact space and $A\subset
C(K)$ is a bounded set of continuous functions over $K$, we define
the pseudometric $d_A$ on $K$ as $d_A(x,y) = \sup_{f\in
A}|f(x)-f(y)|$. If $X$ is a topological space, $d:X\times
X\To\mathbb{R}$ is a map, and $\Delta$ is a positive real number,
it is said that $d$ $\Delta$-\emph{fragments} $X$ if for each
subset $L$ of $X$ there is a relative open subset $U$ of $L$ of
$d$-diameter less than or equal to
$\Delta$.\\

\begin{defn}
A compact space $K$ is said to be countably lower fragmentable if there are bounded subsets
$\{A_{n,p} : n,p\in\mathbb{N}\}$ of $C(K)$ such that $C(K) =
\bigcup_{n\in\mathbb{N}}A_{n,p}$ for every $p\in\mathbb{N}$, and the pseudometric
$d_{A_{n,p}}$ $\frac{1}{p}$-fragments $K$.\\
\end{defn}

This is the definition as it appears in~\cite{FabHeiMat}. However,
variable $p$ is superfluous in it. If the sets $A_{n,1}$ exist, it
is sufficient to define $A_{n,p} = \{\frac{1}{p}f : f\in
A_{n,1}\}$.\\

On the other hand, we recall a concept introduced by
Namioka~\cite{Namioka}: For a topological space $K$, a set
$L\subset K\times K$ is said to be an \emph{almost neighborhood of
the diagonal} if it contains the diagonal $\Delta_K = \{(x,x) :
x\in K\}$ and satisfies that for every nonempty subset $X$ of $K$
there is a nonempty relative open subset $U$ of $X$ such that
$U\times
U\subset L$. The use of this was suggested to us by I. Namioka and simplifies our original proof.\\

\begin{thm}\label{clf equals qRN}
For a compact subset $K$ of $[0,1]^\Gamma$ the following are
equivalent:
\begin{enumerate}
\item $K$ is quasi Radon-Nikodým compact \item $K$ is countably
lower fragmentable. \item There are subsets $\Gamma_{n,p}$,
$n,p\in\mathbb{N}$, of $\Gamma$ such that $d_{\Gamma_{n,p}}$
$\frac{1}{p}$-fragments $K$ for every $n,p\in\mathbb{N}$.\\
\end{enumerate}
\end{thm}

PROOF: Suppose $K$ is quasi Radon-Nikodým compact and let $\phi$
be a lower semicontinuous quasi metric which fragments $K$. Then,
we just define
$$A_{n,p}=\left\{f\in C(K) : |f(x)-f(y)|<\frac{1}{p}
\text{ whenever } \phi(x,y)\leq\frac{1}{n}\right\}\cap\left\{f :
\|f\|_{\infty}\leq n\right\}$$ Clearly, $d_{A_{n,p}}$
$\frac{1}{p}$-fragments $K$ because any subset of $K$ of
$\phi$-diameter less than $\frac{1}{n}$ has $d_{A_{n,p}}$-diameter
less than $\frac{1}{p}$, and we know that $\phi$ fragments $K$. On
the other hand, for a fixed $p\in\mathbb{N}$, in order to prove
that $C(K) = \bigcup_{n\in\mathbb{N}} A_{n,p}$, observe that, if
$f\in C(K)$, then
$$ C_n = \left\{(x,y)\in K\times K : |f(x)-f(y)|\geq \frac{1}{p}\text{ and }
\phi(x,y)\leq\frac{1}{n} \right\}$$ is a decreasing sequence of
compact subsets of $K\times K$ with empty intersection so there is
some $n>\|f\|_{\infty}$ such that $C_n$ is empty, and then, $f\in
A_{n,p}$.\\

That (2) implies (3) is evident, just to take $\Gamma_{n,p} =
A_{n,p}\cap \Gamma$ whenever $A_{n,p}$, $n,p\in\mathbb{N}$ are the
sets in the definition of countably lower fragmentability.\\

Now, suppose (3). For every $n,p\in\mathbb{N}$, since
$d_{A_{n,p}}$ $\frac{1}{p}$-fragments $K$, this means that the set
$C_{n,p} = \{(x,y)\in K\times K :
d_{\Gamma_{n,p}}(x,y)\leq\frac{1}{p}\}$ is an almost neighborhood
of the diagonal which, in addition, is closed. On the other hand,
observe that, for each $n,p\in\mathbb{N}$, $(x,y)\in C_{n,p}$ if
and only if $|x_\gamma-y_\gamma|\leq\frac{1}{p}$ for all
$\gamma\in\Gamma_{n,p}$ so that

$$\bigcap_{n,p\in\mathbb{N}}C_{n,p} = \bigcap_{p\in\mathbb{N}}
\left\{(x,y) : |x_\gamma - y_\gamma|\leq\frac{1}{p}\ \forall
\gamma\in\bigcup_{n\in\mathbb{N}}\Gamma_{n,p}=\Gamma\right\} =
\Delta_K$$

Now, $K$ is quasi Radon-Nikodým by virtue of~\cite[Theorem
1]{Namiokanote}, which states that $K$ is quasi Radon-Nikodým
compact if and only if there is a countable family of closed
almost neighborhoods of the diagonal whose intersection is the
diagonal $\Delta_K$.\hfill$\square$\\

\end{document}